\documentclass[12pt]{amsart}

\usepackage{amssymb}
\newcommand{\K}{\Bbbk}

\newtheorem{lemma}{Lemma}
\newtheorem{theorem}[lemma]{Theorem}

\newtheorem{proposition}[lemma]{Proposition}

\newtheorem{conjecture}[lemma]{Conjecture}

\theoremstyle{definition}
\newtheorem{definition}[lemma]{Definition}

\newtheorem{remark}[lemma]{Remark}
\theoremstyle{remark}
\renewcommand{\theequation}%
{\arabic{section}.\arabic{lemma}.\arabic{equation}}

\begin{document}

\title{Some cases of the Eisenbud-Green-Harris conjecture}

\author{Giulio Caviglia}
\address{Department of Mathematics, U.C. Berkeley, Berkeley, CA 94720}
\email{caviglia@math.berkeley.edu}

\author{Diane Maclagan}
\address{Department of Mathematics, Rutgers University, Hill Center - Busch Campus, Piscataway, NJ 08854}
\email{maclagan@math.rutgers.edu}

\begin{abstract}The Eisenbud-Green-Harris conjecture states that a
homogeneous ideal in $\K[x_1,\dots,x_n]$ containing a homogeneous
regular sequence $f_1,\dots,f_n$ with $\deg(f_i)=a_i$ has the same
Hilbert function as an ideal containing $x_i^{a_i}$ for $1 \leq i \leq
n$.  In this paper we prove the Eisenbud-Green-Harris conjecture when
$a_j > \sum_{i=1}^{j-1} (a_i-1)$ for all $j>1$.  This result was independently obtained by the two authors.
\end{abstract}

\maketitle

\section{Introduction}

Let $\mathbf{a}=(a_1,\dots,a_n)\in \mathbb N^n$, where $2 \leq a_1
\leq a_2 \leq \dots \leq a_n$.  The
following conjecture was originally made by Eisenbud, Green, and
Harris~\cite{EGH} in the case where all the $a_i$ are $2$:

\begin{conjecture}[$\mathrm{EGH}_{\mathbf{a},n}$]  \label{c:egh} Let $S=\K[x_1,\dots,x_n]$, where $\K$ is a field, and let $I$ be a
homogeneous ideal in $S$ containing a regular sequence $f_1,\dots,f_n$
of degrees $\deg(f_i)=a_i$.  Then $I$ has the same Hilbert function as
an ideal containing $\{ x_i^{a_i} : 1 \leq i \leq n \}$.
\end{conjecture}

In this paper we prove the following theorem.

\begin{theorem}
The conjecture $\mathrm{EGH}_{\mathbf{a},n}$ is true if each $a_j$ for
$2 \leq j \leq n$ is larger than $\sum_{i=1}^{j-1} (a_i-1)$.
\end{theorem}

 En-route to proving the theorem we show
the equivalence of the Conjecture $\mathrm{EGH}_{\mathbf{a},n}$ with
similar conjectures for regular sequences of length less than $n$.
Firstly, if $I$ contains a regular sequence of degrees $a_1\leq
\dots \leq a_r$ for $r<n$, and $\mathrm{EGH}_{\mathbf{a},n}$ holds for
all $\mathbf{a}'$ with $\mathbf{a}'_i=a_i$ for $1 \leq i \leq r$ then
$I$ has the same Hilbert function as an ideal containing $x_i^{a_i}$
for $1 \leq i \leq r$.  More importantly,

\begin{proposition} \label{p:prop9}
If $\mathrm{EGH}_{\mathbf{a},n}$ holds then every homogeneous ideal in
the polynomial ring $\K[x_1,\dots,x_m]$ with $m>n$ that contains a
regular sequence $f_1,\dots,f_n$ with $\deg(f_i)=a_i$ has the same
Hilbert function as an ideal containing $x_i^{a_i}$ for $1 \leq i \leq
n$.
\end{proposition}

One of the original motivations of the Conjecture
$\mathrm{EGH}_{\mathbf{a},n}$ was to refine the bounds given by
Macaulay on the possible values of the Hilbert function $H(I,d+1)$ of
an ideal $I$ when $H(I,d)$ is specified in the case where $I$ contains
a regular sequence of degrees $\mathbf{a}$.  A consequence of
Proposition~\ref{p:prop9} is that any known case of
$\mathrm{EGH}_{\mathbf{a},n}$ leads to a refined bound.  For example,
since $\mathrm{EGH}_{\mathbf{a},2}$ holds for any $\mathbf{a}$, the
knowledge that an ideal $I$ in a polynomial ring $\K[x_1,\dots,x_n]$
contains a regular sequence of length two gives a smaller bound on
$H(I,d+1)$ given $H(I,d)$ than that given by Macaulay, since the
Hilbert function of such an ideal must agree with that of an ideal
containing $\{x_1^{a_1}, x_2^{a_2}\}$.

\vspace{0.5cm}

{\bf Acknowledgements:} The results in this paper were independently
discovered by the two authors.  The first author thanks David
Eisenbud, and the second author thanks Greg Smith, both for many
useful conversations.  Maclagan was partially supported by NSF grant
DMS-0500386.

\section{Background}

We first recall the definition of a lexicographic ideal, which will
play a key role in our proof.

\begin{definition}
  An ideal $I \subseteq S$ is {\em lexicographic} with respect to $x_1
  > x_2 > \dots > x_n$ if whenever $x^v \succ_{lex} x^u$, with $x^u
  \in I$, and $\deg(x^u)=\deg(x^v)$, then $x^v \in I$.  There is a
  unique lexicographic ideal with a given Hilbert function (see, for
  example, Section 4.2 of \cite{BrunsHerzog}).
\end{definition}

We will need the following generalization of the notion of a
lexicographic ideal, which is due to Clements and Lindstr\"om
\cite{ClementsLindstrom} (see also \cite{GreeneKleitman}).  We use
here the notation $\mathbb N^n_{\leq}$ for the set of all sequences
$(a_1,\dots,a_n) \in \mathbb N^n$ with $a_1 \leq a_2 \leq \dots \leq
a_n$.  We denote such a sequence by $\mathbf{a}$.

\begin{definition}
Let $\succ_{lex}$ be the lexicographic order with $x_1>x_2> \dots >
  x_n$.  An ideal $I\subseteq S$ is lex-plus-powers with respect to
  the sequence $\mathbf{a} \in \mathbb N^n_{\leq}$ if $\langle x_1^{a_1},
  \dots, x_n^{a_n} \rangle \subseteq I$, and if $x^u \in I$, with $u_i
  < a_i$ for $1 \leq i \leq n$, then for every $x^v$ with
  $\deg(x^u)=\deg(x^v)$, $v_i < a_i$ for $1 \leq i \leq n$ and $x^v
  \succ_{lex} x^u$ we have $x^v \in I$.
\end{definition}

An ideal $I$ is lex-plus-powers if we can write $I=\langle
x_1^{a_1},\dots,x_n^{a_n} \rangle + J$, where $J$ is a lexicographic
ideal.  Note that the order on the variables in the lexicographic
order is forced by the ordering of $\mathbf{a}$.  Indeed, there may be
no lex-plus-powers ideal with respect to another order of the
variables.  A simple example is given by the sequence
$\mathbf{a}=(2,3)$, and the ideal $I=\langle x^2, xy, y^4 \rangle$,
which is lex-plus-powers for the order $x \succ y$.  There is no ideal
of the form $\langle x^2, y^3\rangle + J$ where $J$ is a lexicographic
ideal with respect to the order $y \succ x$ with the same Hilbert
function as $I$.  We also note that this definition differs slightly
from that in the work of Richert~\cite{Richert} and
Francisco~\cite{Francisco}, as we do not require the monomials
$x_i^{a_i}$ to be {\em minimal} generators of the lex-plus-powers
ideal.  See also \cite{EvansCharalambous}, \cite{Cooper},
\cite{PeevaMermin1}, \cite{PeevaMermin2}, and
\cite{PeevaMerminStillman}.

Clements and Lindstr\"om \cite{ClementsLindstrom} show that for any
homogeneous ideal containing $\langle x_1^{a_1},\dots, x_n^{a_n}
\rangle$ there is a lex-plus-powers ideal for the sequence
$\mathbf{a}$ with the same Hilbert function.  Thus
$\mathrm{EGH}_{\mathbf{a},n}$ may be restated in an equivalent form
as: if $I$ contains a regular sequence of degrees $\mathbf{a}$ then
there is lex-plus-powers ideal with respect to $\mathbf{a}$ with the
same Hilbert function.

Note that by taking $a_{k+1},\dots, a_n$ to be arbitrarily large, we
actually get the fact that if $I$ contains $\langle x_1^{a_1},\dots,
x_k^{a_k} \rangle$ for $k<n$, then there is a lex-plus-powers ideal
containing $\langle x_1^{a_1}, \dots, x_k^{a_k} \rangle$ with the
same Hilbert function. Moreover, since the Hilbert function of a
monomial ideal is independent of the base field, by flat extension
of $\K$, we can assume without loss of generality that $\vert \K
\vert= \infty.$

For convenience we list the following well-known facts about regular
sequences that we will use.

\begin{lemma}
\label{l:regseqfacts}
\begin{enumerate}
\item If $J=\langle f_1,\dots,f_r \rangle \subseteq \K[x_1,\dots,x_n]$
is generated by a regular sequence with $\deg(f_i)=a_i$, then $J$ has the same Hilbert function as $\langle x_1^{a_1},\dots, x_r^{a_r} \rangle$.

\item A regular sequence of homogeneous elements in $\K[x_1,\dots,x_n]$
remains a regular sequence after any permutation.

\item If $f_1,\dots,f_r$ is a regular sequence in $\K[x_1,\dots,x_n]$,
with $\K$ infinite, and $a_{r+1},\dots, a_m >0$ for $m \leq n$, then
there exist $f_{r+1},\dots,f_m$ such that $\deg(f_i)=a_i$, and
$f_1,\dots, f_m$ is a regular sequence.

\end{enumerate}

\end{lemma}

The following proposition, which is  Theorem 3 of \cite{DGO}  (see also Corollary
5.2.19 in \cite{Migliore}), will also be useful.

\begin{proposition} \label{p:liaisonhilbfcn}
Let $J=\langle f_1,\dots,f_n \rangle$ be an ideal in $S$ generated
by a regular sequence with $\deg(f_i)=a_i$.  Let $I$ be an ideal
containing $J$ and let $s=\sum_{i=1}^n (a_i-1).$  Then
$$H(S/J,t)=H(S/I,t)+H(S/(J:I),s-t)$$ for $0 \leq
t \leq s $.
\end{proposition}

\section{Proof of the main theorem}
 Our approach to the main theorem will involve the following
relaxation of the conjecture, where we do not assume that the
regular sequence is full length.

\begin{conjecture}
[$\mathrm{EGH}_{n,\mathbf{a},r}$] Let $S=\K[x_1,\dots,x_n]$, and let $I$ be a homogeneous ideal in $S$
containing a regular sequence $f_1,\dots,f_r$ of degrees
$\deg(f_i)=a_i$.  Then $I$ has the same Hilbert function as an ideal
containing $\{ x_i^{a_i} : 1 \leq i \leq r \}$.
\end{conjecture}

 Note that for $\mathrm{EGH}_{n,\mathbf{a},r}$, $\mathbf{a}$ lies
in $\mathbb N^r_{\leq}$.  In this notation
$\mathrm{EGH}_{\mathbf{a},n}$ is $\mathrm{EGH}_{n,\mathbf{a},n}$. As
the following two propositions show, $\mathrm{EGH}_{n,\mathbf{a},r}$
is not actually a generalization of the original conjecture
$\mathrm{EGH}_{\mathbf{a},n}$.

\begin{proposition}
  Fix $n>0$ and $\mathbf{a} \in \mathbb N^r_{\leq}$.  If
  $\mathrm{EGH}_{\mathbf{a}',n}$ holds for all $\mathbf{a}' \in
  \mathbb N^n_{\leq}$ with $a'_i=a_i$ for $1 \leq i \leq r$, then
  $\mathrm{EGH}_{n,\mathbf{a},r}$ holds.
\end{proposition}

\begin{proof}
Suppose that $I \subseteq S$ contains a regular sequence
$f_1,\dots,f_r$ with $\deg(f_i)=a_i$.  Fix $a_{r+1} \leq a_{r+2} \leq
\dots \leq a_n$ with $a_{r+1}>a_r$, and find $f^1_{r+1},\dots f^1_n$
with $\deg(f^1_i)=a_i$ and $f_1,\dots, f_r, f^1_{r+1},\dots,f^1_n$ a
regular sequence.  This is possible by Lemma~\ref{l:regseqfacts}.  Let
$I_1=I+\langle f_{r+1}^1,\dots f_n^1 \rangle$.  Note that
$H(S/I_1,k)=H(S/I,k)$ for $1 \leq k < a_{r+1}$.  Since
$\mathrm{EGH}_{n}$ holds, $I_1$ has the same Hilbert function as an
ideal containing $x_i^{a_i}$ for $1 \leq i \leq n$.  Let $I_{lex}^1$
be the lex-plus-powers ideal with respect to $x_1 > x_2 > \dots > x_n$
with the same Hilbert function as $I_1$ and let $K_1$ be the ideal
generated by those monomials in $I_{lex}^1$ of degree less than
$a_{r+1}$.  Then $K_1$ also has the same Hilbert function as $I$ in
degrees less than $a_{r+1}$, and contains $x_1^{a_1},\dots,x_r^{a_r}$.

Now replace $f^1_{r+1},\dots, f^1_n$ by $f^2_{r+1},\dots,f^2_n$, with
$\deg(f^2_i)=2a_i$, and $f_1,\dots,f_r,f^2_{r+1},\dots,f^2_n$ a
regular sequence.  Set $I_2=I+\langle f^2_{r+1},\dots,f^2_n \rangle$,
let $I_{lex}^2$ be the lex-plus-powers ideal with respect to $x_1 >x_2
> \dots > x_n$ containing
$\{x_1^{a_1},\dots,x_r^{a_r},x_{r+1}^{2a_{r+1}},\dots, x_n^{2a_n} \}$
with the same Hilbert function as $I_2$, and let $K_2$ be the ideal
generated by those monomials in $I_{lex}^2$ of degree less than
$2a_{r+1}$.  Note that $K_2$ has the same Hilbert function as $I$ in
degrees less than $2a_{r+1}$, and contains $x_1^{a_1}, \dots,
x_r^{a_r}$.  Also $K_1$ and $K_2$ agree in degrees less than
$a_{r+1}$, since in degree $k$ their standard monomials are the
$H(S/I,k)$ lexicographically smallest monomials that are not divisible
by $x_1^{a_1},\dots,x_r^{a_r}$.  Thus $K_1 \subseteq K_2$.

We can continue in this manner, choosing $f^j_{r+1},\dots,f^j_n$ with
$\deg(f_i)=ja_i$ completing a regular sequence.  In this manner we get
an increasing sequence $K_1\subseteq K_2 \subseteq K_3 \subseteq
\dots$ of monomial ideals.  Since $S$ is noetherian, this sequence
must eventually stabilize, so there is some $N$ with $K_j$ equal to
$K_N$ for all $j \geq N$.  By construction $K_N$ has the same Hilbert
function as $I$ and contains $x_1^{a_1},\dots, x_r^{a_r}$, so is the
desired ideal.
\end{proof}

\begin{proposition} \label{p:rimpliesnr} Let $\mathbf{a} \in \mathbb
  N^r_{\leq}$.  If $\mathrm{EGH}_{\mathbf{a},r}$ holds for some $r$, then
  $\mathrm{EGH}_{n,\mathbf{a},r}$ holds for all $n\geq r$.
\end{proposition}

\begin{proof}
  It suffices to show that if $\mathrm{EGH}_{n-1,\mathbf{a},r}$ holds
  for some $n-1 \geq r$, then $\mathrm{EGH}_{n,\mathbf{a},r}$ holds.

Suppose that $\mathrm{EGH}_{n-1,\mathbf{a},r}$ holds, and let $f_1,\dots,f_r$ be
a regular sequence contained in an ideal $I \subseteq
S=\K[x_1,\dots,x_{n-1},y]$.  We need to show that there is an ideal
$K\subseteq S$ containing $x_1^{a_1},\dots,x_r^{a_r}$ with the same
Hilbert function as $I$.

By Lemma~\ref{l:regseqfacts} we can find a linear form $g$ such
that $f_1,\dots,f_r,g$, and thus $g,f_1,\dots,f_r$, are regular
sequences.  Let $N>0$ be such that $(I:g^{\infty})=(I:g^N)$.

Note that $R=S/\langle g \rangle$ is isomorphic to a polynomial ring
in $n-1$ variables, and $f_1,\dots,f_r$ descend to a regular sequence
$\overline{f}_1,\dots, \overline{f}_r$ in $R$.  We will construct the
desired ideal $K$ by slices.  Let $I_0=I+\langle g \rangle$, and for
$1 \leq j \leq N$ let $I_j=(I:g^j)+\langle g \rangle$.  Then for $0
\leq j \leq N$ the ideal $I_j$ regarded as an ideal of $R$ contains
$\overline{f}_1,\dots,\overline{f}_r$, so by the induction hypothesis there is an ideal in
$\K[x_1,\dots,x_{n-1}]$ containing $x_1^{a_1},\dots, x_r^{a_r}$ with
the same Hilbert function as $I_j$.  Let $M_j$ be the lex-plus-powers
ideal in $\K[x_1,\dots,x_{n-1}]$ containing
$x_1^{a_1},\dots,x_r^{a_r}$ with this Hilbert function.  Let $K_j
\subseteq S$ be the set of monomials $\{ x^uy^j : x^u \in M_j \}$.
Let $K_{\infty}=\{ x^uy^{N+j} : j \geq 1, x^u \in M_N \}$.

Let $K$ be the ideal generated by the monomials in $K_0,\dots, K_N$.
We claim that $K$ has the desired Hilbert function, and contains
$x_i^{a_i}$ for $1 \leq i \leq r$.  The latter claim is immediate,
since $x_i^{a_i} \in K_0$ for $ 1 \leq i \leq r$.  We prove the former
claim by first showing that $\{ x^uy^j : x^uy^j \in K \}
=\cup_{j=0}^{N} K_j \cup K_{\infty}$.  Note that the sets $K_0,\dots,
K_N,K_{\infty}$ are pairwise disjoint.

To see this we show that the right-hand set is closed under taking
factors, so is the set of monomials in a monomial ideal.  It thus
suffices to show that if $x^uy^j$ is not in the right-hand set, then
neither is any monomial of the form $x^uy^{j-1}$ or $x^uy^j/x_i$.  The
latter is immediate from the definition of $K_j$, as if $x^uy^j \not
\in K_j$, with $j \leq N$, then $x^u \not \in M_j$, so $x^u/x_i \not
\in M_j$, and so $x^uy^j/x_i \not \in K_j$.  If $j>N$, then $x^uy^j
\not \in K_{\infty}$ means that $x^u \not \in M_N$, and so $x^u/x_i
\not \in M_N$, and thus $x^uy^j/x_i \not \in K_{\infty}$.  To see the
former claim, we note that since $(I:g^{j-1}) \subseteq (I:g^{j})$, we
have $I_{j-1} \subseteq I_{j}$, and thus $M_{j-1} \subseteq M_{j}$.
So if $j\leq N$ and $x^uy^j \not \in K_j$, then $x^u \not \in M_j$,
and thus $x^u \not \in M_{j-1}$, so $x^uy^{j-1} \not \in K_{j-1}$.  If
$j>N$, then $x^uy^j \not \in K_{\infty}$ means that $x^u \not \in
K_N$, so $x^u \not \in M_N$, and thus $x^uy^{j-1} \not \in K_N \cup
K_{\infty}$.  This shows that the right-hand side set is the set of
monomials in a monomial ideal, and since $K$ is by definition the
monomial ideal generated by these monomials, we have the equality.

We finish by showing that $K$ has the correct Hilbert function.
Recall that $(I:g^j)=(I:g^N)$ for $j \geq N$, so $(I:g^j)+\langle g
\rangle = I_N$ for such $j$.  Note that the number of monomials
$x^uy^j$ not in $K_j$ (or $K_{\infty}$ if $j>N$) of degree $t$ is the
number of monomials $x^u$ of degree $t-j$ not in $M_j$ (or $M_N$), so
we have the following formula for the Hilbert function of $S/K$:
\begin{eqnarray*}
H(S/K,t)&=&\sum_{j=0}^N H(S/M_j,t-j) + \sum_{j=N+1}^t H(S/M_N,t-j)
    \nonumber \\
&=&\sum_{j=0}^N H(S/I_j,t-j) + \sum_{j=N+1}^t H(S/I_N,t-j)
    \nonumber \\
&=& \sum_{j=0}^t H(S/((I:g^j),g),t-j). \\
\end{eqnarray*}

By considering the short exact sequence
$$0 \rightarrow S/(J:g) \longrightarrow S/J \longrightarrow S/(J,g)
\rightarrow 0,$$ we see that for any ideal $J$ we have
$H(S/J,t)=H(S/(J:g),t-1)+H(S/(J,g),t)$.  Applying this repeatedly we
see that $H(S/I,t)=\sum_{j=0}^t H(S/((I:g^j),g),t-j)=H(S/K,t)$,
completing the proof.
\end{proof}
Conjecture~\ref{c:egh} can be thought as a conjecture on the growth of
ideals containing a regular sequence of given degrees. More precisely
let $I$ be a homogeneous ideal containing a regular sequence of
degrees given by $\mathbf{a}$ and let $d$ be a non-negative
integer. Note that there exists a unique lex-plus-power ideal $J$ of
the form $J=\langle x_1^{a_1}, \dots, x_n^{a_n} \rangle$+$L$, where
$L$ is a lexicographic ideal generated by monomials of the same degree
$d$, with the property that $H(I,d)=H(J,d)$.  It is known that
$\mathrm{EGH}_{n,\mathbf{a}}$ holds if and only if for any $I,d$ and
$J$, defined as above, the condition $H(I,d+1)\geq H(J,d+1)$ is also
satisfied.  We therefore specialize the conjecture at any single
degree in the following way.

\begin{definition} Let $d$ be a non-negative integer. We say that $\mathrm{EGH}_{\mathbf{a},n}(d)$
holds if for any homogeneous ideal $I\subseteq \K[x_1,\dots,x_n]$
containing a regular sequence of degrees  $\mathbf{a} \in \mathbb
{N}^n_{\leq}$ there exists a homogeneous ideal $J$ containing  $\{
x_i^{a_i} : 1 \leq i \leq r \}$ such that $H(I,d)=H(J,d)$ and
$H(I,d+1)=H(J,d+1).$
\end{definition}
Note that by the Clements-Lindstr\"om
Theorem~\cite{ClementsLindstrom} we can assume that $J$ is a
lex-plus-powers ideal with respect to $\mathbf{a}$.
Thus $\mathrm{EGH}_{n,\mathbf{a}}$ holds if and only if
$\mathrm{EGH}_{\mathbf{a},n}(d)$  holds for all non-negative
integers $d.$  The following Lemma shows  a symmetry in
Conjecture~\ref{c:egh}.

\begin{lemma}\label{l:symmetry} Let $\mathbf{a} \in \mathbb
{N}^n_{\leq}$ and let $s= \sum_{i=1}^n (a_i-1).$ Then for a
non-negative integer $d$ we have that
$\mathrm{EGH}_{\mathbf{a},n}(d)$ holds if and only if
$\mathrm{EGH}_{\mathbf{a},n}(s-d-1)$ holds.
\end{lemma}

\begin{proof} Assume that $\mathrm{EGH}_{\mathbf{a},n}(d)$  holds.
Let $I$ be a homogeneous ideal of $S=\K[x_1,\dots,x_n]$ containing a
regular sequence $f_1,\dots,f_n$ of degrees $\mathbf{a} \in
\mathbb{N}^n_{\leq}.$ Let $F=\langle f_1,\dots, f_n \rangle$ and let
$I_1=(F:I)$. By $\mathrm{EGH}_{\mathbf{a},n}(d)$ there exists an ideal
$J$ containing $M=\langle x_1^{a_1},\dots, x_n ^{a_n} \rangle$ such
that $H(I_1,d)=H(J,d)$ and $H(I_1,d+1)=H(J,d+1)$.  Set $J_1$ equal to
$(M:J)$.  Note $J_1$ contains $M$ so by Lemma~\ref{l:regseqfacts} and
Proposition~\ref{p:liaisonhilbfcn} we have $H(I,s-d)=H(J_1,s-d)$ and
$H(I,s-d-1)=H(J_1,s-d-1)$.
\end{proof}

\begin{theorem} \label{t:mainthm}
The conjecture $\mathrm{EGH}_{\mathbf{a},n}$ is true if $a_j >
\sum_{i=1}^{j-1} (a_j-1)$ for all $j>1$.
\end{theorem}

\begin{proof}The proof is by induction on $n$.  When $n=1$ the hypothesis
$\mathrm{EGH}_{a,1}$ states that if a homogeneous ideal $I \subset
\K[x_1]$ contains an element $f$ of degree $a_1$ then it has the same
Hilbert function as an ideal containing $x_1^{a_1}$.  This is
immediate, since $x_1^{a_1}$ is the only homogeneous polynomial of
degree $a_1$ in $\K[x_1]$.  We now assume that
$\mathrm{EGH}_{\mathbf{a}',n-1}$ holds, where $\mathbf{a}' \in \mathbb
N^{n-1}_{\leq}$ is the projection of $\mathbf{a}$ onto the first $n-1$
coordinates. Let $s=\sum_{i=1}^n (a_i-1).$ By Lemma~\ref{l:symmetry}
it is enough to prove $\mathrm{EGH}_{\mathbf{a},n}(d)$ for $0 \leq d
\leq \lfloor (s-1)/2 \rfloor.$ Let $I\subseteq \K[x_1,\dots,x_n]$ be a
homogeneous ideal containing a regular sequence of degrees
$\mathbf{a}.$ By Proposition~\ref{p:rimpliesnr} we have that
$\mathrm{EGH}_{n,\mathbf{a'},n-1}$ holds, so there is an ideal J,
containing $x_1^{a_1},\dots,x_{n-1}^{a_{n-1}}$ with the same Hilbert
function as $I$. Since, by assumption, $a_n> \lfloor (s-1)/2
\rfloor+1$ we deduce that $H(I,d)=H(J+(x_n^{a_n}),d)$ for all $0\leq d
\leq \lfloor (s-1)/2 \rfloor +1.$
\end{proof}

\begin{remark}
Note that Theorem~\ref{t:mainthm} implies
$\mathrm{EGH}_{\mathbf{a},2}$ for any $\mathbf{a}=(a_1,a_2)$, which was
already known.
\end{remark}

\end{document}